\documentclass[12pt]{amsart}
\usepackage{amsfonts,amsmath,amsthm,oldgerm,amssymb,amscd}

\usepackage[all,cmtip]{xy}

\usepackage{fancyhdr}
\usepackage{psfrag}
\usepackage{graphicx}

\usepackage[indent]{caption}
\usepackage{indentfirst}
\usepackage{float}
\usepackage{latexsym}
\usepackage{graphics}
\usepackage{verbatim}

\newcommand{\Q}{{\mathbb Q}}
\newcommand{\R}{{\mathbb R}}
\newcommand{\Z}{{\mathbb Z}}
\newcommand{\C}{{\mathbb C}}

\newcommand{\G}{{\mathbb G}}

\newcommand{\fg}{{\mathfrak g}}
\newcommand{\fk}{{\mathfrak k}}
\newcommand{\ft}{{\mathfrak t}}
\newcommand{\fb}{{\mathfrak b}}
\renewcommand{\O}{{\mathcal O}}

\newtheorem{theorem}{Theorem}

\newtheorem*{thm}{Theorem}

\theoremstyle{definition}

\numberwithin{question}{section}

\theoremstyle{remark}

\begin{document}

\title[Convexity of tensor products]{On a convexity property of tensor products of irreducible, rational representations of $SL(n)$ }

\author{Hariharan Narayanan}
\address{School of Technology and Computer Science, Tata Institute of Fundamental
  Research, Homi Bhabha Road, Mumbai, India 400005.\\
email: hariharan.narayanan@tifr.res.in}  

\author{C. S. Rajan}
\address{School of Mathematics, Tata Institute of Fundamental
  Research,  Homi Bhabha Road, Mumbai, India 400005.\\
email: rajan@math.tifr.res.in}

\thanks{The authors acknowledge the support of the Department of Atomic Energy, Government of India under project no. 12-RD-TFR-RT14001. HN Acknowledges support from a Ramanujan fellowship and a Swarna Jayanti  fellowship instituted by the government of India.}
\keywords{} 
\subjclass{Primary 20G05; Secondary 22E45, 05E10}

\begin{abstract} 
The aim of this note is to point out a convexity property with respect to the root lattice for the support of the highest weights that occur in a tensor product of irreducible rational representations of $SL(n)$ (or more generally any Lie group over which the saturation property holds) over the complex numbers. The observation is a consequence of the convexity properties of the saturation cone and the validity of the saturation conjecture for $SL(n)$. 
\end{abstract}

\maketitle

\section{Main Theorem} 
Let $G$ be a connected reductive algebraic group over $\C$. Choose a Borel subgroup $B$ of $G$ and a maximal torus $T$ of $B$. Denote by $\Lambda=X^*(T)=\mbox{Hom}(T, \G_m)$ (resp.  $\Lambda^\vee=X_*(T)=\mbox{Hom}(\G_m, T)$ the character (resp. the cocharacter) group of $T$. Let $\langle.,.\rangle$ denote the natural pairing $X^*(T)\times X_*(T)\to \Z$. Denote by  $\Phi\subset \Lambda$ (resp. $\Phi^\vee\subset \Lambda^\vee$) the set of roots (resp. coroots) of $(G,T)$, and $\Phi^+$ the subset of positive roots formed by the roots of $(B,T)$. Let $Q$ denote the root lattice, the subgroup of $\Lambda$ generated by $\Phi$. 

The dominant chamber $\Lambda^+$ of $\Lambda$ is the collection of 
characters $\lambda\in \Lambda$ such that $\langle\lambda, \alpha^\vee\rangle$ are non-negative for $\alpha\in \Phi^+$, where for a root $\alpha$, $\alpha^\vee\in \Phi^\vee$ denotes the co-root corresponding to $\alpha$. The irreducible, rational, linear  representations are parametrized by their highest weights. For an element $\lambda\in \Lambda^+$, denote by $V(\lambda)$ the corresponding irreducible representation with highest weight $\lambda$. 
We recall that a Lie group $G$ is said to have the {\it saturation property} if the following is true.
Whenever $\lambda_1, \cdots, \lambda_s$ are highest weights of $G$ such that 
 $\lambda_1+\cdots+\lambda_s$\ belongs to the root lattice. Assume further that for some natural number $m$, the space of $G$-invariants for the tensor product representation $V(m\lambda_1)\otimes \cdots \otimes  V(m\lambda_s)$ is non-zero.  Then, 
 \[ [V(\lambda_1)\otimes \cdots \otimes  V(\lambda_s)]^G\neq 0.\]

The main aim of this note is to observe the following convexity property of irreducible, rational representations of any Lie group that satisfies the saturation property, in particular, $SL(n)$, over complex numbers:
\begin{theorem}\label{thm:main}
Let $V_1, \cdots, V_r$ be  irreducible rational representations of 
$SL(n)$ over complex numbers. Suppose there exists a natural number  $k$,  dominant weights $\nu$ and $\nu+k\gamma$ belonging to $\Lambda^+$ with $\gamma$ an element in the root lattice $Q$,  such that 
the corresponding irreducible representations $V(\nu)$ and $V(\nu+k\gamma)$ occur in the tensor product $V_1\otimes \cdots \otimes V_r$. 

Then for any natural number $l$ with $0\leq l\leq k$, the irreducible representation $V(\nu+l\gamma)$ with highest weight 
$\nu+l\gamma$ occurs in $V_1\otimes \cdots \otimes V_r$.
\end{theorem}

For any natural number $l$ with $0\leq l\leq k$, the weight 
$\nu+l\gamma= \frac{(k-l)\nu}{k}+\frac{l(\nu+k\gamma)}{k}$ is a rational, convex linear combination of the characters  $\nu$ and $\nu+k\gamma$. In particular, the characters  $\nu+l\gamma$ are all dominant. 

We recall, that the irreducible, rational representations of $SL(n)$  are indexed by integer partitions. The correspondence is given as follows: suppose   $\lambda=m_1\omega_1+\cdots +m_{n-1}\omega_{n-1}$ is a highest weight, where $\omega_i$ are the fundamental representations of $SL(n)$ corresponding respectively to the $i$-th exterior power representatation of the standard representation and  $m_i$ non-negative integers. The corresponding partition is given by $a_1\geq a_2\geq \cdots a_{n-1}\geq a_n=0$, where $a_1=m_1+\cdots+m_{n-1}, ~a_2=m_2+\cdots+m_{n-1}, \cdots , ~a_{n-1}=m_{n-1}$ and $a_n=0$. 
The root vectors of $SL(n)$ are vectors in $\R^n$ with exactly two nonzero entries, one being $1$ and the other being $-1$. The root lattice $Q$ is the subgroup of $\Z^n$, consisting of those integral elements in $\Z^n$, whose sum of its co-ordinates is zero. Theorem~\ref{thm:main}  may be rephrased more explicitly as follows:

\begin{thm}
Let $\lambda^{(1)}, \dots, \lambda^{(r)}$ be a collection of $r$ integer partitions, with $\lambda^{(i)} = (a^{(i)}_1, \dots a^{(i)}_{n})$ being a partition with non-increasing nonnnegative parts. Let $V_1, \dots, V_r$ be the corresponding irreducible representations of $SL(n)$.
Suppose there exists a natural number  $k$, non-decreasing integral partitions $\nu$ and $\nu+k\gamma$ with $\gamma$ belonging to the  root lattice $Q$  such that the corresponding irreducible representations $V(\nu)$ and $V(\nu+k\gamma)$ occur in the tensor product $V_1\otimes \cdots \otimes V_r$ of irreducible finite dimensional representations of $SL(n)$. 

Then for any natural number $l$ with $0\leq l\leq k$, the irreducible representation $V(\nu+l\gamma)$ with highest weight corresponding to the partition
$\nu+l\gamma$ occurs in $V_1\otimes \cdots \otimes V_r$.
\end{thm}

The proof of Theorem~\ref{thm:main} follows from  (a) the convexity properties of the moment map and its relationship to the restrictions of characters of compact groups following the work of Heckman (\cite{He}), Guillemin-Sternberg (\cite{GS}), Mumford-Ness (\cite{N}) Kirwan (\cite{K}),  and (b) the validity of the saturation conjecture for $SL(n)$ due to Knutson-Tao (\cite{KT}). We recall these results following the expositions (\cite{Br, SK2}). 

In what follows, for an algebraic group $M$ over $\C$,  we also use the  symbol $M$ to denote the complex points of $M$. Let $K$ be a maximal compact subgroup of $G$, such that $T_K=T\cap K$ is a maximal torus of $K$.  We denote respectively the Lie algebras of $K$ and $T_K$ by $\fk$ and  $\ft_K$. Since $T_K$ is Zariski dense in $T$, the characters $\lambda\in \Lambda$ are determined by their restrictions to $T_K$. We continue to denote by $\lambda$ the differential of the restriction of $\lambda$ to $T_K$: these belong to the dual space $\ft_K^*\subset \fk^*$ of $\ft_K$. Let $W$ be the Weyl group and let $\ft_K^+$ denote the fundamental chamber for the action of $W$, given by the choice of $B$. 

Given elements $\lambda_1, \cdots, \lambda_r$ belonging to $\fk^*$, consider the coadjoint orbit $$\O=K^r(\lambda_1, \cdots, \lambda_r)\subset (\fk^*)^r$$ of the compact group $K^r$. The space $\O$ is a compact, connected manifold, equipped with a natural symplectic structure for which the moment map is given by the inclusion of $\O$ in $(\fk^*)^r$.

Restricting the action of $K^r$ to the diagonal action of $K$ on $\O$, the moment map $p$ for the action is given by composing the inclusion with the sum $(\fk^*)^r\to \fk^*$, sending $(X_1, \cdots, X_r)\in  (\fk^*)^r$ to $X_1+\cdots +X_r$. By the convexity theorem (\cite[Theorem 1.3]{Br}), the image $p(\O)$ intersected with the dominant chamber is a convex polytope. 

Further, in the case that $\lambda_i$ arise from $X^*(T)\otimes \Q$ and are dominant, the equivalence of symplectic reduction with that of the GIT quotient together with applying the Borel-Weil theorem, imply that 
the set of dominant $\nu\in X^*(T)\otimes \Q$ that lie in the image $p(\O)$ are precisely those for which there exists a natural number $m$ such that the weights $m\nu, m\lambda_1,\cdots, m\lambda_r$ are dominant characters (belonging to $\Lambda^+$), such that the corresponding highest weight representation $V(m\nu)$ occurs in the tensor product $V(m\lambda_1)\otimes \cdots \otimes  V(m\lambda_r)$ (\cite[Theorem 1.3]{Br}). 

In the hypothesis of Theorem \ref{thm:main}, let $\lambda_1, \cdots,  \lambda_r$ be respectively the highest weights of the representations $V_1, \cdots, V_r$ of $SL(n,\C)$. Given the above discussion, the dominant weights $\nu$ and $\nu+k\gamma$ belong to the image $p(\O)$. Since the image $p(\O)$ is convex, the dominant weights $\nu+l\gamma$ also belong to 
$p(\O)$ for $0\leq l\leq k$. In particular, given any $l$ with  $0\leq l\leq k$, there exists a natural number $m$ such that the representation $V(m(\nu+l\gamma))$ occurs in $V(m\lambda_1)\otimes \cdots \otimes  V(m\lambda_r)$. 

Let $w_0$ be the longest element in the Weyl group, sending $\Phi^+$ to $-\Phi^+$. The highest weight of the dual representation $V(\lambda)^*$ is given by $-w_0\lambda$. Hence, we get that the space of $G=SL(n)$ invariants of the tensor product 
$$V(m\lambda_1)\otimes \cdots \otimes  V(m\lambda_r)\otimes V(-mw_0(\nu+l\gamma))$$
is non-zero. 

We recall now the saturation theorem for $SL(n)$ proved by Knutson and Tao (\cite{KT}, \cite[Theorem 13]{SK2}):
\begin{theorem}\label{thm:sat}
Suppose $\lambda_1, \cdots, \lambda_s$ are highest weights of $SL(n)$ such that 
 $\lambda_1+\cdots+\lambda_s$\ belongs to the root lattice. Assume further that for some natural number $m$, the space of $G$-invariants for the tensor product representation $V(m\lambda_1)\otimes \cdots \otimes  V(m\lambda_s)$ is non-zero.  Then, 
 \[ [V(\lambda_1)\otimes \cdots \otimes  V(\lambda_s)]^G\neq 0.\]

\end{theorem}

The weights that occur in a highest weight representation $V(\lambda)$ are of the form $\lambda-\beta$ for some $\beta$ belonging to the (positive) root lattice. Hence the weights that occur in the tensor product $V(\lambda_1)\otimes \cdots \otimes  V(\lambda_r)\otimes V(-w_0\nu)$ are of the form $\lambda_1+\cdots+\lambda_r-w_o\nu-\beta$ for some  $\beta\in Q^+$. Since this tensor product has a non-zero space of $SL(n)$-invariants, the $0$-weight occurs in the tensor product. In particular, $\lambda_1+\cdots+\lambda_r-w_0\nu$ belongs to the root lattice, and hence for any integer $l$, the weight $\lambda_1+\cdots+\lambda_r-w_0\nu+l\gamma)$ belongs to the root lattice. 

Hence, Theorem \ref{thm:main} follows 
by appealing to the saturation theorem. 

\section{Remarks}

1. The weights that occur in the restriction of a highest weight representation $V(\lambda)$ of $G$ to a maximal torus $T$,  is the intersection of the translate $\lambda+Q$ of the root lattice $Q$ by $\lambda$,  with that of the convex span of the extremal weights $W\lambda$ (\cite[Chapter VIII, Section 7]{Bou}). In particular, the convexity property with respect to the root lattice holds for the above restriction. 

Suppose $\mu$ is a weight that occurs in the restriction of $V(\lambda)$. Consider the $\alpha$-chain through $\mu$: let $p$ (resp. $q$) be the largest non-negative integer such that $\mu+p\alpha$ (resp. $\mu-q\alpha$) is a weight of $V(\lambda)$. For $\alpha$ a root, it is known that $q-p=\mu(\alpha^\vee)$. 
These results follow from the representation theory of the Lie algebra $sl(2)$ associated to the root $\alpha$. 

It would be interesting to know whether similar methods or results hold in the context of Theorem \ref{thm:main}.

2. Saturation conjecture need not hold for groups other than $SL(n)$. It is conjectured by Kapovich and Millson, that it should hold for a tensor product $V(\lambda)\otimes V(\mu)$ provided the weights $\lambda$ and $\mu$ are regular (\cite[Appendix, Conjecture 85]{SK2}).

More generally, one can consider the general restriction problem,  to study the restriction 
of representations of a given reductive group $\hat{G}$ to a reductive subgroup $G$ (\cite{Br}). 

Using the explicit description of the classical branching laws, as given in (\cite[Chapter 8.3]{GW}), it can be seen that the restriction of an irreducible representation of $GL(n)$ to $GL(n-1)$ (where $GL(n-1)$ embeds in the usual way) 
(and for such other classical pairs) satisfies the convexity property with respect to the root lattice as given above. In this case, much more is known: the weights of the restriction satisfy a log-concavity property (\cite{Ok}). However, log-concavity fails in general, for example for tensor products of irreducible representations of $GL(n)$ (\cite{CDW}).  It will be interesting to know whether it will hold provided the highest weights are regular. 
It was recently shown (see \cite{HMMS}) that for Kostka numbers, which are the weight multiplicities for irreducible representations of $GL(n)$, log-concavity does hold.

3. Suppose $V$ is a vector space of dimension $2l+1\geq 5$ and $\Psi$ be a non-degenerate quadratic form on $V$. Let $G=O(V, \Psi)$ be the orthogonal group of the quadratic form. The group $G$ is of type $B_l$. Consider the standard representation $V$ of $G$, with highest weight $\omega_1$ in the notation of \cite[Ch. VIII, Section 13.2]{Bou}. The tensor product $V\otimes V$ breaks up as $S^2V\oplus \Lambda^2V$. The  exterior square representation $\Lambda^2V$ is irreducible with highest weight $\omega_2$. The symmetric square representation $S^2V$ breaks up as a sum of an irreducible representation $V(2\omega_1)$ with highest weight $2\omega_1$ and a copy of a trivial representation corresponding to the quadratic form. Further, the weight $\omega_1$ is a also a root. 

This example shows that the convexity property  with respect to the restriction of the root lattice fails. It also gives a counterexample to the saturation property. Thus it seems that the convexity property of the support with respect to the root lattice of the 
ambient group is linked to the validity of saturation type conjectures, and we possibly cannot expect a direct proof of the convexity property,  avoiding the use of the validity of saturation conjecture for $SL(n)$. 

4. Another question that naturally comes up: if the collection of weights that occur in the restriction of a tensor product as above, is indeed convex with respect to the root lattice of the ambient group, then what are its
extreme points? The Horn conjectures give only the codimension one subspaces, the boundary 
faces. 

Given highest weights $\lambda$ and $\mu$ for $SL(n)$, the validity of the PRV conjecture implies that for any element $w$ of the Weyl group of $SL(n)$, the irreducible representation having $\lambda+w\mu$ as an extreme weight occurs in the tensor product $V(\lambda)\otimes V(\mu)$ (\cite{SK, MPR}).  Some of the PRV components and its generalizations give extreme points. It is known that for large $n$, the PRV type analogues that occur in the Horn inequalities  do not suffice to determine the facets that describe the Horn polytope.  It will be interesting to know the convex span of the PRV components with respect to the root lattice.

Morever, the fact that convex polytopes corresponding to Kostka numbers, (i.e., the polytopes whose integer points correspond to semistandard Young Tableaux with weight $\lambda$ and content $\mu$) have vertices that are non-integral \cite{DLM}, (though rational) suggests that the problem of characterizing all vertices of these polytopes  is difficult.

\noindent
{\bf Acknowledgement.} We thank Rekha Biswal, Michel Brion, Apoorva Khare, Joshua Kiers, Dipendra Prasad  and B. Ravinder for useful discussions. The second author thanks MPIM,  Bonn for a visit during May 2019, for an excellent working environment during which part of this work was carried out.

\end{document}